\documentclass[10pt,leqno]{article}
\usepackage{amsfonts,amsmath}
\newtheorem{theorem}{Theorem}[section]
\newtheorem{definition}{Definition}[section]
\newtheorem{proposition}{Proposition}[section]

\begin{document}
\title{On the $k$-Semispray of Nonlinear Connections in $k$-Tangent Bundle Geometry}
\author{Florian Munteanu}
\date{Department of Applied Mathematics,\\
University of Craiova, Romania\\
munteanufm@gmail.com}

\maketitle

\begin{abstract}
In this paper we present a method by which is obtained a sequence of $k$-semisprays and two
sequences of nonlinear connections on the $k$-tangent bundle
$T^kM$, starting from a given one. Interesting particular cases appear for Lagrange and Finsler spaces of order $k$.
\end{abstract}

\textbf{AMS} \textbf{Subject Classification:} 53C05, 53C60.

\textbf{Key words:} $k$-tangent bundle, $k$-semispray, nonlinear connection,
Lagrange space of order $k$, Finsler space of order $k$.
\section{Introduction}
\par \noindent Classical Mechanics have been entirely geometrized in terms of
symplectic geometry and in this approach there exists certain
dynamical vector field on the tangent bundle $TM$ of a manifold
$M$ whose integral curves are the solutions of the Euler-Lagrange
equations. This vector field is usually called \textit{spray} or
\textit{second-order differential equation (SODE)}. Sometimes it is called \textit{semispray} and the
term \textit{spray} is reserved to homogeneous  second-order
differential equations (\cite{miron2}, \cite{roman-roy}). Let us remember that a
SODE on $TM$ is a vector field on $TM$ such that $JC=C$, where $J$
is the almost tangent structure and $C$ is the canonical Liouville field (\cite{leon_rodrigues}, \cite{miron1}).

In \cite{grif1}, \cite{grif2}, \cite{grif3} J. Grifone studies the relationship among SODEs, nonlinear connections and the autonomous Lagrangian formalism. In paper \cite{muntgh} Gh. Munteanu and Gh. Piti\c{s} also studied the relation between sprays and nonlinear connectiosn on $TM$.   This study was extended to the non-autonomous case by M. de Le\'{o}n and P. Rodgrigues (\cite{leon_rodrigues}). Also, important results for singular non-autonomous case was obtained in  \cite{fmlp}. In this paper, following the ideas of papers \cite{muntfl4}, \cite{muntfl5}, \cite{muntgh} and \cite{fmlp} we will extend the study of the  relationship  between sprays and nonlinear connections to the $k$-tangent bundle of a manifold $M$. The study of the geometry of this  $k$-tangent bundle was by introduced by R. Miron (\cite{miron2}, \cite{miron2p}, \cite{miron3}). For this case the $k$-spray represent a system of ordinary differential equations of $k+1$ order.

\section{The $k$-Semispray of a Nonlinear Connection}

Let $M$ be a real $n$-dimensional manifold of class $C^{\infty}$
and $(T^kM,\pi ^k,M)$ the bundle of accelerations of order $k$. It
can be identified with the $k$-osculator bundle or $k$-tangent
bundle (\cite{miron2}, \cite{miron3}).

A point $u\in T^kM$ will be written by
$u=(x,y^{(1)},...,y^{(k)})$, $\pi ^k(u)=x$, $ x\in M$. The
canonical coordinates of $u$ are $(x^i,y^{(1)i},...,y^{(k)i})$,
$i=\overline{1,n}$, where $y^{(1)i}=\displaystyle\frac
1{1!}\displaystyle\frac{dx^i}{dt}$, ..., $
y^{(2)i}=\displaystyle\frac
1{k!}\displaystyle\frac{d^kx^i}{dt^k}$. A transformation of local
coordinates $(x^i,y^{(1)i},...,y^{(k)i})\rightarrow
(\widetilde{x}^i,\widetilde{y} ^{(1)i},...,\widetilde{y}^{(k)i})$
on $(k+1)n$-dimensional manifold $T^kM$ is given by
\begin{equation}
\left\{
\begin{tabular}{l}
$\widetilde{x}^i=\widetilde{x}^i(x^1,...,x^n)$, $rang$ $\left(
\displaystyle\frac{
\partial \widetilde{x}^i}{\partial x^j}\right) =n$, \\
$\widetilde{y}^{(1)i}=\displaystyle\frac{\partial
\widetilde{x}^i}{\partial x^j}y^{(1)j}$,
\\
$2\widetilde{y}^{(2)i}=\displaystyle\frac{\partial
\widetilde{y}^{(1)i}}{\partial x^j}
y^{(1)j}+2\displaystyle\frac{\partial
\widetilde{y}^{(1)i}}{\partial y^{(1)j}}y^{(2)j}$,
\\
$........................................................................
$ \\
$k\widetilde{y}^{(k)i}=\displaystyle\frac{\partial
\widetilde{y}^{(k-1)i}}{\partial x^j}
y^{(1)j}+2\displaystyle\frac{\partial
\widetilde{y}^{(k-1)i}}{\partial y^{(1)j}} y^{(2)j}+\cdots
+k\displaystyle\frac{\partial \widetilde{y}^{(k-1)i}}{\partial
y^{(k-1)j}} y^{(k)j}.$
\end{tabular}
\right.  \label{mfs1}
\end{equation}
A local coordinates change (\ref{mfs1}) transforms the natural basis\\
$\left\{ \displaystyle\frac \partial {\partial
x^i},\displaystyle\frac \partial {\partial y^{(1)i}},\cdots
,\displaystyle\frac \partial {\partial y^{(k)i}}\right\}_u$ of the
tangent space $T_uT^kM$ by the rule:
\begin{equation}
\left\{
\begin{array}{ccccccccc}
\displaystyle\frac \partial {\partial x^i} & = &
\displaystyle\frac{\partial \widetilde{x}^j}{\partial
x^i}\displaystyle\frac
\partial {\partial \widetilde{x}^j} & + & \displaystyle\frac{\partial
\widetilde{y}^{(1)j}}{\partial x^i}\displaystyle\frac \partial
{\partial \widetilde{y} ^{(1)j}} & + & \cdots & + &
\displaystyle\frac{\partial \widetilde{y}^{(k)j}}{\partial x^i
}\displaystyle\frac \partial {\partial \widetilde{y}^{(k)j}}, \\
\displaystyle\frac \partial {\partial y^{(1)i}} & = &  &  &
\displaystyle\frac{\partial \widetilde{y} ^{(1)j}}{\partial
y^{(1)i}}\displaystyle\frac
\partial {\partial \widetilde{y}^{(1)j}} & + & \cdots & + &
\displaystyle\frac{\partial \widetilde{y}^{(k)j}}{\partial
y^{(1)i}}\displaystyle\frac
\partial {\partial \widetilde{y}^{(k)j}}, \\
\vdots &  &  &  &  &  &  &  &  \\
\displaystyle\frac \partial {\partial y^{(k)i}} & = &  &  &  &  &
&  & \displaystyle\frac{\partial \widetilde{y}^{(k)j}}{\partial
y^{(k)i}}\displaystyle\frac
\partial {\partial \widetilde{y }^{(k)j}} .
\end{array}
\right.  \label{mfs2}
\end{equation}
The distribution $V_1:u\in T^kM\rightarrow V_{1,u}\subset T_uT^kM$
generated by the tangent vectors $\left\{ \displaystyle\frac
\partial {\partial y^{(1)i}},\cdots ,\displaystyle\frac \partial {\partial
y^{(k)i}}\right\} _u$ is a vertical distribution on the bundle
$T^kM$. Its local dimension is $kn$. Similarly, the distribution
$V_2:u\in T^kM\rightarrow V_{2,u}\subset T_uT^kM$ generated by
$\left\{ \displaystyle\frac \partial {\partial y^{(2)i}},\cdots
,\displaystyle\frac \partial {\partial y^{(k)i}}\right\}_u$ is a
subdistribution of $V_1$ of local dimension $(k-1)n$. So, by this
procedure on obtains a sequence of integrable distributions $V_1
\supset V_2 \supset \cdots \supset V_{k}$. The last distribution
$V_k$ is generated by $\left\{ \displaystyle\frac
\partial {\partial y^{(k)i}}\right\}_u$ and $\dim V_k=n$ (\cite{miron2}).

Hereafter, we consider the open submanifold
\[
\widetilde{T^kM}=T^kM\setminus \{\mathbf{0}\}=\left\{
(x,y^{(1)},...,y^{(k)})\in T^kM|rank\ ||y^{(1)i}||=1\right\} \, ,
\]
where $\mathbf{0}$ is the null section of the projection
$\pi^k:T^kM\rightarrow M$.

The following operators in algebra of functions
$\mathcal{F}(T^kM)$
\begin{equation}
\begin{array}{l}
\overset{1}{\Gamma }=y^{(1)i}\displaystyle\frac \partial {\partial y^{(k)i}}, \\
\overset{2}{\Gamma }=y^{(1)i}\displaystyle\frac \partial {\partial
y^{(k-1)i}}+2y^{(2)i}\displaystyle\frac \partial {\partial y^{(k)i}}, \\
.................................................................... \\
\overset{k}{\Gamma }=y^{(1)i}\displaystyle\frac \partial {\partial
y^{(1)i}}+2y^{(2)i}\displaystyle\frac \partial {\partial
y^{(2)i}}+\cdots +ky^{(k)i}\displaystyle\frac
\partial {\partial y^{(k)i}}
\end{array}
\label{mfs3}
\end{equation}
are $k$ vector fields, globally defined on $T^kM$ and linearly
independent on the manifold $\widetilde{T^kM}=T^kM\setminus
\{\mathbf{0}\}$, $ \overset{1}{\Gamma}$ belongs of distribution
$V_k$, $\overset{2}{\Gamma}$ belongs of distribution $V_{k-1}$,
..., $\overset{k}{\Gamma}$ belongs of distribution $V_1$ (see
\cite{miron2}). $\overset{1}{\Gamma }$ , $\overset{2}{\Gamma }$,
..., $\overset{k}{\Gamma }$ are called \textit{Liouville vector
fields}.

In applications we shall use also the following nonlinear
operator, which is not a vector field,
\begin{equation}
\Gamma =y^{(1)i}\displaystyle\frac \partial {\partial
x^i}+2y^{(2)i}\displaystyle\frac \partial {\partial
y^{(1)i}}+\cdots +ky^{(k)i}\displaystyle\frac \partial {\partial
y^{(k-1)i}}. \label{mfs4}
\end{equation}
Under a coordinates transformation (\ref{mfs1}) on $T^kM$,
$\Gamma$ changes as follows:
\begin{equation}
\Gamma =\widetilde{\Gamma }+\left\{
y^{(1)i}\displaystyle\frac{\partial \widetilde{y}
^{(k)j}}{\partial x^i}+\cdots
+ky^{(k)i}\displaystyle\frac{\partial \widetilde{y}^{(k)j}}{
\partial y^{(k-1)i}}\right\} \displaystyle\frac \partial {\partial \widetilde{y}^{(k)j}}.
\label{mfs5}
\end{equation}
A $k$-\textit{tangent structure} $J$ on $T^kM$ is defined as
usually (\cite{miron2}) by the following
$\mathcal{F}(T^kM)$-linear mapping $J:\mathcal{X}(T^kM)\rightarrow
\mathcal{X} (T^kM)$:
\begin{equation}
\begin{array}{c}
J\left( \displaystyle\frac \partial {\partial x^i}\right)
=\displaystyle\frac
\partial {\partial y^{(1)i}},J\left( \displaystyle\frac \partial {\partial
y^{(1)i}}\right) =\displaystyle\frac \partial
{\partial y^{(2)i}},..., \\
J\left( \displaystyle\frac \partial {\partial y^{(k-1)i}}\right)
=\displaystyle\frac \partial {\partial y^{(k)i}},J\left(
\displaystyle\frac
\partial {\partial y^{(k)i}}\right) =0.
\end{array}
\label{mfs6}
\end{equation}
$J$ is a tensor field of type $(1,1)$, globally defined on $T^kM$.
\begin{definition}
(\cite{miron2}) \textit{A} $k$\textit{-semispray} on $T^kM$ is a
vector field $S \in \mathcal{X}(T^kM)$ with the property
\begin{equation}
JS=\overset{k}{\Gamma }.  \label{mfs7}
\end{equation}
\end{definition}
Obviously, there not always exists a $k$-semispray, globally
defined on $T^kM$. Therefore the notion of local $k$-semispray is
necessary. For example, if $M$ is a paracompact manifold then on
$T^kM$ there exists local $k$-semisprays (\cite{miron2}).
\begin{theorem}
(\cite{miron2}) i) A $k$-semispray $S$ can be uniquely written in
local coordinates in the form:
\begin{equation}
\begin{array}{lll}
S & = & y^{(1)i}\displaystyle\frac \partial {\partial
x^i}+2y^{(2)i}\displaystyle\frac
\partial {\partial y^{(1)i}}+\cdots +ky^{(k)i}\displaystyle\frac \partial
{\partial y^{(k-1)i}}-
\\
& - & (k+1)G^i(x,y^{(1)},...,y^{(k)})\displaystyle\frac \partial
{\partial y^{(k)i}} \, .
\end{array}
\label{mfs8}
\end{equation}
ii) With respect to (\ref{mfs1}) the coefficients
$G^i(x,y^{(1)},...,y^{(k)})$ change as follows:
\begin{equation}
\begin{array}{lll}
(k+1)\widetilde{G}^i & = & (k+1)G^j\displaystyle\frac{\partial
\widetilde{x}^i}{\partial
x^j}- \\
& - & \left( y^{(1)j}\displaystyle\frac{\partial
\widetilde{y}^{(k)i}}{\partial x^j} +\cdots
+ky^{(k)j}\displaystyle\frac{\partial
\widetilde{y}^{(k)i}}{\partial y^{(k-1)j}} \right) \, .
\end{array}
\label{mfs9}
\end{equation}
iii) If the functions $G^i(x,y^{(1)},...,y^{(k)})$ are given on
every domain of local chart of $T^kM$, so that (\ref{mfs9}) holds,
then the vector field $S$ from (\ref{mfs8}) is a $k$-semispray.
\end{theorem}
Let us consider a curve $c:I\rightarrow M$, represented in a local
chart $(U,\varphi)$ by $x^i=x^i(t)$, $t\in I$. Thus, the mapping
$\widetilde{c} : I \rightarrow T^kM$, given on $(\pi ^k)^{-1}(U)$,
by
\begin{equation}
x^i=x^i(t),y^{(1)i}(t)=\displaystyle\frac
1{1!}\displaystyle\frac{dx^i}{dt}(t),\ldots
,y^{(k)i}(t)=\displaystyle\frac
1{k!}\displaystyle\frac{d^kx^i}{dt^k}(t),t\in I \label{mfs10}
\end{equation}
is a curve in $T^kM$, called the $k$-\textit{extension} to $T^kM$
of the curve $c$.

A curve $c:I\rightarrow M$ is called $k$-\textit{path} of a
$k$-semispray $S$ (from (\ref{mfs8})) if its $k$-extension
$\widetilde{c}$ is an integral curve for $S$, that is
\begin{equation}
\left\{ \displaystyle\frac{dx^i}{dt}=y^{(1)i},\
\displaystyle\frac{dy^{(1)i}}{dt}=2y^{(2)i},\ ...,
\displaystyle\frac{dy^{(k-1)i}}{dt}=ky^{(k)i},\displaystyle\frac{dy^{(k)i}}{dt}
=-(k+1)G^i \, . \right.  \label{mfs11}
\end{equation}
\begin{definition}
The $k$-semispray $S$ is called $k$-\textit{spray} if the
functions $\left( G^i(x,y^{(1)},...,y^{(k)})\right)$ are
$(k+1)$-homogeneous, that is
\[
G^i(x,\lambda y^{(1)},...,\lambda ^ky^{(k)})=\lambda
^{k+1}G^i(x,y^{(1)},...,y^{(k)}),\ \forall \lambda >0.
\]
\end{definition}
Like in the case of tangent bundle, an Euler Theorem holds. That
is, a function $f\in \mathcal{F}(\widetilde{T^kM})$ is
$r$-homogeneous if and only if
\[
\mathcal{L}_{\overset{k}{\Gamma }}f=rf.
\]
Then a $k$-semispray $S$ is a $k$-spray if and only if
\begin{equation}
y^{(1)h}\displaystyle\frac{\partial G^i}{\partial
y^{(1)h}}+2y^{(2)h}\displaystyle\frac{\partial G^i}{
\partial y^{(2)h}}+\cdots +ky^{(k)h}\displaystyle\frac{\partial G^i}{\partial y^{(k)h}}
=(k+1)G^i.  \label{mfs12}
\end{equation}
\begin{definition}
A vector subbundle $NT^kM$ of the tangent bundle $(TT^kM,d\pi
^k,M)$ which is supplementary to the vertical subbundle $V_1T^kM$,
\begin{equation}
TT^kM=NT^kM\oplus V_1T^kM  \label{mfs13}
\end{equation}
is called \textit{a nonlinear connection} on $T^kM$.
\end{definition}
The fibres of $NT^kM$ determine a horizontal distribution $N:u\in
T^kM\rightarrow N_uT^kM\subset T_uT^kM$ supplementary to the
vertical distribution $V_1$, that is
\begin{equation}
T_uT^kM=N_uT^kM\oplus V_{1,u}T^kM,\ \forall u\in T^kM.  \label{mfs13p}
\end{equation}
The dimension of horizontal distribution $N$ is $n$.

If the base manifold $M$ is paracompact then on $T^kM$ there
exists the nonlinear connections (\cite{miron2}).

There exists a unique local basis, adapted to the horizontal
distribution $N$, $\left\{ \displaystyle\frac \delta {\delta
x^i}\right\}_{i=\overline{1,n}}$, such that $d\pi ^k\left(
\displaystyle\frac \delta {\delta x^i}|_u\right)
=\displaystyle\frac
\partial {\partial x^i}|_{\pi ^k(u)}$, $i=1$, ..., $n$. More over,
on each domain of local chart of $T^kM$ there exists the functions
$\underset{(1)}{N_j^i}$, $\underset{(2)}{N_j^i}$, ...,
$\underset{(k)}{N_j^i}$ such that
\begin{equation}
\displaystyle\frac \delta {\delta x^i}=\displaystyle\frac \partial
{\partial x^i}-\underset{(1)}{ N_i^j}\displaystyle\frac \partial
{\partial y^{(1)j}}-\cdots -\underset{(k)}{N_i^j}
\displaystyle\frac \partial {\partial y^{(k)j}} \, .
\label{mfs14}
\end{equation}
The functions $\underset{(1)}{N_j^i}$, $\underset{(2)}{N_j^i}$,
...,$ \underset{(k)}{N_j^i}$ are called \textit{the primal
coefficients} of the nonlinear connection $N$ and under a
coordinates transformation (\ref{mfs1}) on $T^kM$ this
coefficients are changing by the rule:
\begin{equation}
\left\{
\begin{array}{lll}
\underset{(1)}{\widetilde{N}_m^i}\displaystyle\frac{\partial
\widetilde{x}^m}{\partial x^j} & = & \displaystyle\frac{\partial
\widetilde{x}^i}{\partial x^m}\underset{(1)}{
N_j^m}-\displaystyle\frac{\partial \widetilde{y}^{(1)i}}{\partial x^j}, \\
\underset{(2)}{\widetilde{N}_m^i}\displaystyle\frac{\partial
\widetilde{x}^m}{\partial x^j} & = & \displaystyle\frac{\partial
\widetilde{x}^i}{\partial x^m}\underset{(2)}{
N_j^m}+\displaystyle\frac{\partial \widetilde{y}^{(1)i}}{\partial
x^m}\underset{(1)}{
N_j^m}-\displaystyle\frac{\partial \widetilde{y}^{(2)i}}{\partial x^j},..., \\
\underset{(k)}{\widetilde{N}_m^i}\displaystyle\frac{\partial
\widetilde{x}^m}{\partial x^j} & = & \displaystyle\frac{\partial
\widetilde{x}^i}{\partial x^m}\underset{(k)}{
N_j^m}+\displaystyle\frac{\partial \widetilde{y}^{(1)i}}{\partial
x^m}\underset{(k-1)}{ N_j^m}+\cdots +\displaystyle\frac{\partial
\widetilde{y}^{(k-1)i}}{\partial x^m}
\underset{(1)}{N_j^m}-\displaystyle\frac{\partial
\widetilde{y}^{(k)i}}{\partial x^j}.
\end{array}
\right.  \label{mfs15}
\end{equation}
Conversely, if on each local chart of $T^kM$ a set of functions
$\underset{(1)}{N_j^i}$, ..., $\underset{(k)}{ N_j^i}$ is given so
that, according to (\ref{mfs1}), the equalities (\ref{mfs15})
hold, then there exists on $T^kM$ a unique nonlinear connection
$N$ which has as coefficients just the given set of function
(\cite{miron2}).

\textit{The local adapted basis} $\left\{ \displaystyle\frac
\delta {\delta x^i},\displaystyle\frac \delta {\delta
y^{(1)i}},\cdots ,\displaystyle\frac \delta {\delta
y^{(k)i}}\right\}_{i=\overline{1,n}}$ is given by (\ref{mfs14})
and
\begin{equation}
\begin{array}{lll}
\displaystyle\frac \delta {\delta y^{(1)i}}=\displaystyle\frac
\partial {\partial y^{(1)i}}- \underset{(1)}{N_i^j}\displaystyle\frac \partial
{\partial y^{(2)j}}-\cdots -
\underset{(k-1)}{N_i^j}\displaystyle\frac \partial {\partial y^{(k)j}},..., &  &  \\
\displaystyle\frac \delta {\delta y^{(k-1)i}}=\displaystyle\frac
\partial {\partial y^{(k-1)i}}- \underset{(1)}{N_i^j}\displaystyle\frac
\partial {\partial y^{(k)j}},\displaystyle\frac \delta {\delta y^{(k)i}}=\displaystyle\frac
\partial {\partial y^{(k)j}} &  &
\end{array}
\label{mfs17}
\end{equation}
and \textit{the dual basis} (or \textit{the adapted cobasis}) of
adapted basis is\\
$\left\{ \delta x^i,\delta y^{(1)i},\ldots
,\delta y^{(k)i}\right\}_{i= \overline{1,n}} \,$, where $\delta
x^i=dx^i$ and
\begin{equation}
\left\{
\begin{array}{lll}
\delta y^{(1)i} & = & dy^{(1)i}+\underset{(1)}{M_j^i}dx^j, \\
\delta y^{(2)i} & = &
dy^{(2)i}+\underset{(1)}{M_j^i}dy^{(1)j}+\underset{
(2)}{M_j^i}dx^j,..., \\
\delta y^{(k)i} & = &
dy^{(k)i}+\underset{(1)}{M_j^i}dy^{(k-1)j}+\cdots +
\underset{(k)}{M_j^i}dx^j
\end{array}
\right.  \label{mfs18}
\end{equation}
and
\begin{equation}
\left\{
\begin{array}{lll}
\underset{(1)}{M_j^i} & = & \underset{(1)}{N_j^i}, \\
\underset{(2)}{M_j^i} & = &
\underset{(2)}{N_j^i}+\underset{(1)}{N_m^i}
\underset{(1)}{M_j^m},\ ..., \\
\underset{(k)}{M_j^i} & = &
\underset{(k)}{N_j^i}+\underset{(k-1)}{
N_m^i}\underset{(1)}{M_j^m}+\cdots
+\underset{(1)}{N_m^i}\underset{ (k-1)}{M_j^m}.
\end{array}
\right.  \label{mfs19}
\end{equation}
Conversely, if the adapted cobasis $\left\{ \delta x^i,\delta
y^{(1)i},\ldots ,\delta y^{(k)i}\right\}_{i=\overline{1,n}}$ is
given in the form (\ref{mfs18}), then the adapted basis $\left\{
\displaystyle\frac \delta {\delta x^i},\displaystyle\frac \delta
{\delta y^{(1)i}},\cdots ,\displaystyle\frac \delta {\delta
y^{(k)i}}\right\}_{i=\overline{1,n}}$ is expressed in the form
(\ref {mfs17}), where
\begin{equation}
\left\{
\begin{array}{lll}
\underset{(1)}{N_j^i} & = & \underset{(1)}{M_j^i}, \\
\underset{(2)}{N_j^i} & = &
\underset{(2)}{M_j^i}-\underset{(1)}{N_m^i}
\underset{(1)}{M_j^m},\ ..., \\
\underset{(k)}{N_j^i} & = &
\underset{(k)}{M_j^i}-\underset{(k-1)}{
N_m^i}\underset{(1)}{M_j^m}-\cdots
-\underset{(1)}{N_m^i}\underset{ (k-1)}{M_j^m}.
\end{array}
\right. \, .  \label{mfs19p}
\end{equation}
The functions $\underset{(1)}{M_j^i}$, $\underset{(2)}{M_j^i}$,
..., $ \underset{(k)}{M_j^i}$ are called \textit{the dual
coefficients} of the nonlinear connection $N$.

A nonlinear connection $N$ is complete determined by a system of
functions $\underset{(1)}{M_j^i}$, ..., $\underset{(k)}{M_j^i}$
which is given on each domain of local chart on $T^kM$, so that,
according to (\ref{mfs1}), the relations hold:
\begin{equation}
\left\{
\begin{array}{lll}
\underset{(1)}{M_j^m}\displaystyle\frac{\partial
\widetilde{x}^i}{\partial x^m} & = & \displaystyle\frac{\partial
\widetilde{x}^m}{\partial x^j}\underset{(1)}{\widetilde{M}
_m^i}+\displaystyle\frac{\partial \widetilde{y}^{(1)i}}{\partial x^j}, \\
\underset{(2)}{M_j^m}\displaystyle\frac{\partial
\widetilde{x}^i}{\partial x^m} & = & \displaystyle\frac{\partial
\widetilde{x}^m}{\partial x^j}\underset{(2)}{\widetilde{M}
_m^i}+\displaystyle\frac{\partial \widetilde{y}^{(1)m}}{\partial
x^j}\underset{(1)}{ \widetilde{M}_m^i}+\displaystyle\frac{\partial
\widetilde{y}^{(2)i}}{\partial x^j},...,
\\
\underset{(k)}{M_j^m}\displaystyle\frac{\partial
\widetilde{x}^i}{\partial x^m} & = & \displaystyle\frac{\partial
\widetilde{x}^m}{\partial x^j}\underset{(k)}{\widetilde{M}
_m^i}+\displaystyle\frac{\partial \widetilde{y}^{(1)m}}{\partial
x^j}\underset{(k-1)}{ \widetilde{M}_m^i}+\cdots
+\displaystyle\frac{\partial \widetilde{y}^{(k-1)m}}{\partial
x^j}\underset{(1)}{\widetilde{M}_m^i}+\displaystyle\frac{\partial
\widetilde{y}^{(k)i}}{
\partial x^j}.
\end{array}
\right.  \label{mfs20}
\end{equation}
Let $c:I\rightarrow M$ be a parametrized curve on the base
manifold $M$, given by $x^i=x^i(t)$, $t\in I$. If we consider its
$k$-extension $\widetilde{c}$ to $T^kM$, then we say that $c$ is
 \textit{an autoparallel curve} for the nonlinear connection $N$ if its $k$-extension
 $\widetilde{c}$ is an horizontal curve, that is
  $\displaystyle\frac{d\widetilde{c}}{dt}$ belongs to the horizontal distribution.

From (\ref{mfs18}) and
\begin{equation}
\displaystyle\frac{d\widetilde{c}}{dt}=\displaystyle\frac{dx^i}{dt}\displaystyle\frac
\delta {\delta x^i}+\displaystyle\frac{ \delta
y^{(1)i}}{dt}\displaystyle\frac \delta {\delta y^{(1)i}}+\cdots
+\displaystyle\frac{\delta y^{(k)i}}{dt}\displaystyle\frac \delta
{\delta y^{(k)i}} \label{mfs21}
\end{equation}
it result that the autoparallels curves of the nonlinear
connection $N$ with the dual coefficients $\underset{(1)}{M_j^i}$,
..., $\underset{(k)}{M_j^i }$ are characterized by the system of
differential equations (\cite {miron2}):
\begin{equation}
\left\{
\begin{array}{l}
y^{(1)i}=\displaystyle\frac{dx^i}{dt},\
y^{(2)i}=\displaystyle\frac
1{2!}\displaystyle\frac{d^2x^i}{dt^2},\ ...,\
y^{(k)i}=\displaystyle\frac 1{k!}\displaystyle\frac{d^kx^i}{dt^k}, \\
\displaystyle\frac{\delta
y^{(1)i}}{dt}=\displaystyle\frac{dy^{(1)i}}{dt}+\underset{(1)}{M_j^i}\displaystyle\frac{
dx^j}{dt}=0, \\
\displaystyle\frac{\delta
y^{(2)i}}{dt}=\displaystyle\frac{dy^{(2)i}}{dt}+\underset{(1)}{M_j^i}\displaystyle\frac{
dy^{(1)j}}{dt}+\underset{(2)}{M_j^i}\displaystyle\frac{dx^j}{dt}=0, \\
..................................................................... \\
\displaystyle\frac{\delta
y^{(k)i}}{dt}=\displaystyle\frac{dy^{(k)i}}{dt}+\underset{(1)}{M_j^i}\displaystyle\frac{
dy^{(k-1)j}}{dt}+\cdots
+\underset{(k)}{M_j^i}\displaystyle\frac{dx^j}{dt}=0.
\end{array}
\right.  \label{mfs22}
\end{equation}
Now, let be $S=\overset{1}{S}$ a $k$-semispray with the
coefficients $\overset{1}{G^i}=G^i(x,y^{(1)},...,y^{(k)})$ like in
(\ref{mfs8}). Then the set of functions
\begin{equation}
\left\{
\begin{array}{l}
\underset{(1)}{M_j^i}=\displaystyle\frac{\partial G^i}{\partial y^{(k)j}}, \\
\underset{(2)}{M_j^i}=\displaystyle\frac 12\left(
S\underset{(1)}{M_j^i}+\underset{
(1)}{M_m^i}\underset{(1)}{M_j^m}\right) , \\
........................................... \\
\underset{(k)}{M_j^i}=\displaystyle\frac 1k\left(
S\underset{(k-1)}{M_j^i}+
\underset{(1)}{M_m^i}\underset{(k-1)}{M_j^m}\right)
\end{array}
\right.  \label{mfs23}
\end{equation}
gives the dual coefficients of a nonlinear connection $N$
determined only by the $k$-semispray $S$ (see the book
\cite{miron2} of Radu Miron).

Other result, obtained by Ioan Buc\u ataru
(\cite{bucataru}), give a second nonlinear connection $N^{*}$ on
$T^kM$ determined only by the $k$-semispray $S$. That is, the
following set of functions
\begin{equation}
\underset{(1)}{M_j^{*i}}=\displaystyle\frac{\partial G^i}{\partial
y^{(k)j}},\ \underset{(2)}{M_j^{*i}}=\displaystyle\frac{\partial
G^i}{\partial y^{(k-1)j}},\ ...,\
\underset{(k)}{M_j^{*i}}=\displaystyle\frac{\partial G^i}{\partial
y^{(1)j}} \label{mfs24}
\end{equation}
is the set of dual coefficients of a nonlinear connection $N^{*}$.

Let us consider the set of functions
$(\overset{2}{G^i}(x,y^{(1)},...,y^{(k)}))$, given on every domain
of local chart by
\begin{equation}
\overset{2}{G^i}=\displaystyle\frac 1{k+1}\overset{k}{\Gamma
}\overset{1}{G^i}=\displaystyle\frac
1{k+1}y^{(1)h}\displaystyle\frac{\partial
\overset{1}{G^i}}{\partial y^{(1)h}}+\displaystyle\frac
2{k+1}y^{(2)h}\displaystyle\frac{\partial
\overset{1}{G^i}}{\partial y^{(2)h}}+\cdots +\displaystyle\frac
k{k+1}y^{(k)h}\displaystyle\frac{\partial
\overset{1}{G^i}}{\partial y^{(k)h}}. \label{mfs25}
\end{equation}
Using (\ref{mfs5}) we obtain that the functions $\overset{2}{G^i}$
verifies (\ref{mfs9}). So, the functions $\overset{2}{G^i}$
represent the coefficients of a $k$-semispray $\overset{2}{S}$,
\begin{equation}
\begin{array}{lll}
\overset{2}{S} & = & y^{(1)i}\displaystyle\frac \partial {\partial
x^i}+2y^{(2)i}\displaystyle\frac
\partial {\partial y^{(1)i}}+\cdots +ky^{(k)i}\displaystyle\frac \partial {\partial
y^{(k-1)i}}- \\
& - &
(k+1)\overset{2}{G^i}(x,y^{(1)},...,y^{(k)})\displaystyle\frac
\partial {\partial y^{(k)i}}\text{.}
\end{array}
\label{mfs26}
\end{equation}
Obviously, there exists two nonlinear connections on $T^kM$, which
 depend only by the $k$-semispray $\overset{2}{S}$:\\
$\overset{2}{N}$ with the dual coefficients
\begin{equation}
\left\{
\begin{array}{l}
\underset{(1)}{\overset{2}{M_j^i}}=\displaystyle\frac{\partial
G^i}{\partial y^{(k)j}},
\\
\underset{(2)}{\overset{2}{M_j^i}}=\displaystyle\frac 12\left(
S\underset{(1)}{
\overset{2}{M_j^i}}+\underset{(1)}{\overset{2}{M_m^i}}\underset{(1)}{
\overset{2}{M_j^m}}\right) , \\
........................................... \\
\underset{(k)}{\overset{2}{M_j^i}}=\displaystyle\frac 1k\left(
S\underset{(k-1)}{
\overset{2}{M_j^i}}+\underset{(1)}{\overset{2}{M_m^i}}\underset{(k-1)}{
\overset{2}{M_j^m}}\right)
\end{array}
\right.  \label{mfs27}
\end{equation}
and $\overset{2}{N^{*}}$ with the dual coefficients
\begin{equation}
\underset{(1)}{\overset{2}{M_j^{*i}}}=\displaystyle\frac{\partial
\overset{2}{G^i}}{
\partial y^{(k)j}},\ \underset{(2)}{\overset{2}{M_j^{*i}}}=\displaystyle\frac{\partial
\overset{2}{G^i}}{\partial y^{(k-1)j}},\ ...,\
\underset{(k)}{\overset{2 }{M_j^{*i}}}=\displaystyle\frac{\partial
G^i}{\partial y^{(1)j}}.  \label{mfs28}
\end{equation}
By this method is obtained a sequence of $k$-semisprays $\left(
\overset{m}{S}\right)_{m\geq 1}$ and two sequence of nonlinear
connections, $\left( \overset{m}{N}\right)_{m\geq 1}$, $\left(
\overset{m}{ N^{*}}\right)_{m\geq 1}$.

From (\ref{mfs11}), (\ref{mfs22}) and (\ref{mfs25}) we have the
following results:
\begin{proposition}
If $c$ is an autoparallel curve for nonlinear connection
$\overset{1}{N^{*}}$, then $c$ is a $k$-path of $k$-semispray
$\overset{2}{S}$.
\end{proposition}
\begin{theorem}
The following assertions are equivalent:

i) the $k$-semispray $\overset{1}{S}$ is a $k$-spray;

ii) the $k$-paths of $\overset{1}{S}$ and $\overset{2}{S}$
coincide.
\end{theorem}
\begin{theorem}
If $\overset{1}{S}$ is a $k$-spray then
$\overset{1}{\underset{(1)}{M_j^i}}$, ...,
$\overset{1}{\underset{(k)}{M_j^i}}$ (or $\overset{1
}{\underset{(1)}{M_j^{*i}}}$, ...,
$\overset{1}{\underset{(k)}{M_j^{*i}}}$) are homogeneous functions
of degree $1$, $2$, ..., $k$, respectively. The same property have
the primal coefficients $\overset{1}{\underset{(1)}{N_j^i}}$, ...,
$\overset{1}{ \underset{(k)}{N_j^i}}$ (or
$\overset{1}{\underset{(1)}{N_j^{*i}}}$, ...,
$\overset{1}{\underset{(k)}{N_j^{*i}}}$).
\end{theorem}
We remark that the converse of this proposition is generally not
valid and we have the result:
\begin{theorem}
If $\overset{1}{S}$ is a $k$-spray then the sequence $\left(
\overset{m}{S}\right)_{m\geq 1}$ is constant and the sequences
$\left( \overset{m}{N}\right)_{m\geq 1}$ , $\left(
\overset{m}{N^{*}}\right)_{m\geq 1}$ are constant.
\end{theorem}
\section{The $k$-Semispray of a Nonlinear Connection in a Lagrange Space of Order $k$}
\textit{A Lagrangian of order} $k$ is a mapping $
L:T^kM\rightarrow \mathbf{R}$. $L$ is called
\textit{differentiable} if it is of $C^\infty$-class on
$\widetilde{T^kM}$ and continuous on the null section of the
projection $\pi^k:T^kM\rightarrow M$.

The Hessian of a differentiable Lagrangian $L$, with respect to
the variables $y^{(k)i}$ on $\widetilde{T^kM}$ is the matrix
$||2g_{ij}||$, where
\begin{equation}
g_{ij}=\displaystyle\frac 12\displaystyle\frac{\partial
^2L}{\partial y^{(k)i}\partial y^{(k)j}}. \label{mfs29}
\end{equation}
We have that $g_{ij}$ is a $d$-tensor field on the manifold
$\widetilde{T^kM}$, covariant of order $2$, symmetric (see
\cite{miron2}).

If
\begin{equation}
rank\ ||g_{ij}||=n,\ \ \mbox{ on } \ \ \widetilde{T^kM}
\label{mfs30}
\end{equation}
we say that $L(x,y^{(1)},...,y^{(k)})$ is a \textit{regular} (or
\textit{nondegenerate}) Lagrangian.

The existence of the regular Lagrangians of order $k$ is proved
for the case of paracompacts manifold $M$ in the book
\cite{miron2} of Radu Miron.
\begin{definition} (\cite{miron2})
We call \textit{a Lagrange space of order} $k$ a pair
$L^{(k)n}=(M,L)$, formed by a real $n$-dimensional manifold $M$
and a regular differentiable Lagrangian of order $k$,
$L:(x,y^{(1)},...,y^{(k)})\in T^kM\rightarrow
L(x,y^{(1)},...,y^{(k)})\in \mathbf{R}$, for which the quadratic
form $\Psi =g_{ij}\xi^i \xi^j$ on $\widetilde{T^kM}$ has a
constant signature.
\end{definition}
$L$ is called \textit{the fundamental function} and $g_{ij}$
\textit{the fundamental} (or \textit{metric}) \textit{tensor
field} of the space $L^{(k)n}$.

It is known that for any regular Lagrangian of order $k$,
$L(x,y^{(1)},...,y^{(k)})$, there exists a $k$-semispray $S_L$
determined only by the Lagrangian $L$ (see \cite{miron2}). The
coefficients of $S_L$ are given by
\begin{equation}
(k+1)G^i=\displaystyle\frac 12g^{ij}\left\{ \Gamma \left(
\displaystyle\frac{\partial L}{\partial y^{(k)j}}\right)
-\displaystyle\frac{\partial L}{\partial y^{(k-1)j}}\right\} \, .
\label{mfs31}
\end{equation}

This $k$-semispray $S_L$ depending only by $L$ will be called
\textit{canonical}. If $L$ is globally defined on $T^kM$, then
$S_L$ has the same property on $\widetilde{T^kM}$.

From (\ref{mfs23}) and (\ref{mfs24}) it result that there exists
two nonlinear connections: \textit{Miron's connection} $N$ and
\textit{Buc\u ataru's connection} $N^{*}$ which depending only by
the Lagrangian $L$. For this reason, both are called
\textit{canonical}.

So, the coefficients of $k$-semisprays $\overset{m}{S}$ and the
coefficients of nonlinear connections $\overset{m}{N}$,
$\overset{m}{N^{*}}$ depend only by the Lagrangian $L$, for any
$m\geq 1$, but their expressions is not attractive for us.

Interesting results appear for Finsler spaces of order $k$.
\begin{definition} (\cite{miron2p})
\textit{A Finsler space of order} $k$ is a pair $F^{\left(
k\right) n}=(M,F)$ formed by a real differentiable manifold $M$ of
dimension $n$ and a function $F:T^kM\rightarrow \mathbf{R}$ having
the following properties:

i) $F$ is differentiable on $\widetilde{T^kM}$ and continuous on
null section $0:M\rightarrow T^kM$;

ii) $F$ is positive;

iii) $F$ is $k$-homogeneous;

iv) the Hessian of $F^2$ with elements
\begin{equation}
g_{ij}=\displaystyle\frac 12\displaystyle\frac{\partial
^2F^2}{\partial y^{\left( k\right) i}\partial y^{\left( k\right)
j}} \label{mfs32}
\end{equation}
is positively defined on $\widetilde{T^kM}$.
\end{definition}
The function $F$ is called \textit{the fundamental function} and
the $d$-tensor field $g_{ij}$ is called \textit{fundamental}(or
\textit{metric}) \textit{tensor field} of the Finsler space of
order $k$, $F^{\left( k\right) n}$.

The class of spaces $F^{\left( k\right) n}$ is a subclass of
spaces $L^{\left( k\right) n}$.

Taking into account the $k$-homogeneity of the fundamental
function $F$ and $2k$-homogeneity of $F^2$ we get:

1. the coefficients $G^i$ of the canonical $k$-semispray
$S_{F^2}$, determined only by the fundamental function $F$,
\begin{equation}
(k+1)G^i=\displaystyle\frac 12g^{ij}\left\{ \Gamma \left(
\displaystyle\frac{\partial F^2}{\partial y^{(k)j}}\right)
-\displaystyle\frac{\partial F^2}{\partial y^{(k-1)j}}\right\} ,
\label{mfs33}
\end{equation}
is $(k+1)$-homogeneous functions, that is $S_{F^2}$ is a
$k$-spray;

2. the dual coefficients of \textit{the Cartan nonlinear
connection} $N$ associated to Finsler space of order $k$,
$F^{(k)n}$ (see \cite{miron2p}),
\begin{equation}
\begin{array}{l}
\underset{(1)}{M_j^i}=\displaystyle\frac{1}{2(k+1)}
\displaystyle\frac
\partial {\partial y^{(k)j}}\left\{ g^{im}\left[ \Gamma \left(
\displaystyle\frac{\partial F^2}{\partial y^{(k)m}}\right) -
\displaystyle\frac{\partial F^2}{\partial y^{(k-1)m}}\right] \right\} , \\
\underset{(2)}{M_j^i}=\displaystyle\frac 12\left(
S_{F^2}\underset{(1)}{M_j^i}+
\underset{(1)}{M_m^i}\underset{(1)}{M_j^m}\right) , \\
............................................. \\
\underset{(k)}{M_j^i}=\displaystyle\frac 1k\left(
S_{F^2}\underset{(k-1)}{M_j^i}+
\underset{(1)}{M_m^i}\underset{(k-1)}{M_j^m}\right) ,
\end{array}
\label{mfs34}
\end{equation}
are homogeneous functions of degree $1$, $2$, ..., $k$,
respectively, and the primal coefficients has the same property;

3. the dual coefficients of Buc\u ataru's connection $N^{*}$
associated to Lagrangian $F^2$ are also homogeneous functions of
degree $1$, $2$, ..., $k$, respectively, and the primal
coefficients has the same property.

Using the previous results, we obtain the results:
\begin{theorem}
If $F^{(k)n}=(M,F)$ is a Finsler space of order $k$, then:

a) the sequence $\left( \overset{m}{S}\right)_{m\geq 1}$ is
constant, $ \overset{1}{S}$ being the canonical $k$-spray
$S_{F^2}$;

b) the sequences of nonlinear connections $\left(
\overset{m}{N}\right)_{m\geq 1}$, $\left(
\overset{m}{N^{*}}\right)_{m\geq 1}$ are constants,
$\overset{1}{N}$ being the Cartan nonlinear connection of
$F^{(k)n}$ and $\overset{1}{N^{*}}$ being the Buc\u ataru's
connection for $L=F^2$.
\end{theorem}
\section{Conclusions}
In this paper was studied the the relation between semisprays and nonlinear connections on the $k$-tangent bundle $T^kM$ of a manifold $M$. This results was generalized by the author from the $2$-tangent bundle $T^2M$ (\cite{muntfl5}). More that, the relationship between SOPDEs and nonlinear connections on the tangent bundle of $k^1$-velocities of a manifold $M$ (i.e. the Whitney sum of $k$ copies of $TM$, $ T^{1}_{k}M = TM \oplus \dots \oplus TM $) was  studied by F. Munteanu in \cite{muntfl2006} (2006)  and by  N. Roman-Roy, M. Salgado, S. Vilarino in \cite{roman-roy} (2011).

\textbf{Acknowledgments.} This research was partially supported by Grant FP7-PEOPLE-2012-IRSES-316338.
%

\end{document}